\documentclass[english]{article}
\usepackage[T1]{fontenc}
\usepackage[latin9]{inputenc}
\usepackage{color,soul}
\usepackage{framed}
\usepackage{amsmath}
\usepackage{amsthm}
\usepackage{amssymb}

\makeatletter

\providecommand{\tabularnewline}{\\}

\numberwithin{equation}{section}
\numberwithin{figure}{section}
\theoremstyle{plain}
\newtheorem{thm}{\protect\theoremname}
  \theoremstyle{remark}
  \newtheorem{claim}[thm]{\protect\claimname}

\usepackage{fullpage}\usepackage{makeidx}\usepackage{amsfonts}\usepackage{latexsym}

\theoremstyle{definition}
\newtheorem{lemma}{Lemma}[section]
\newtheorem{theorem}[lemma]{Theorem}\newtheorem{definition}[lemma]{Definition}\newtheorem{remark}[lemma]{Remark}\usepackage{times}

\usepackage{babel}

\usepackage[blocks]{authblk}
\title{The VOAs generated by two Ising vectors $e$ and $f$ with $\langle e, f\rangle=\frac{1}{2^8}$ or $\frac{3}{2^9}$}

\author{Wen Zheng}
\affil{Department of Mathematics, University of
California, Santa Cruz, CA 95064 USA}

\usepackage{babel}

\usepackage{babel}

\usepackage{babel}

\usepackage{babel}
\providecommand{\claimname}{Claim}
\providecommand{\theoremname}{Theorem}

\usepackage{babel}
\providecommand{\claimname}{Claim}
\providecommand{\theoremname}{Theorem}

\numberwithin{equation}{section}

\usepackage{babel}
\providecommand{\claimname}{Claim}
\providecommand{\theoremname}{Theorem}

\makeatother
\usepackage{babel}
  \providecommand{\claimname}{Claim}
\providecommand{\theoremname}{Theorem}
\def\ha{\frac{1}{2}}
\begin{document}
\definecolor{shadecolor}{rgb}{0,1,0}
\maketitle
\begin{abstract}
 In this paper we study the VOAs generated by two Ising vectors whose inner product is $\frac{1}{2^8}$ or $\frac{3}{2^9}$ and determine they both 
have unique VOA structures.

\end{abstract}
\section{Introduction}
\sethlcolor{green}

In the study of the moonshine vertex operator algebra $V^{\natural}$ constructed in \cite{FLM}, Dong, et al. in \cite{DMZ} showed that
$V^{\natural}$ contains $48$ Virasoro vectors, each Virasoro vector generates a Virasoro vertex operator algebra isomorphic to $L(\ha,0)$ in $V^{\natural}$ and $L(\ha,0)^{\otimes 48}$ is a conformal subalgebra of $V^{\natural}.$ Such a Virasoro vector is called an Ising vector. Later, Miyamoto in \cite{M} constructed a $\tau$-involution $\tau_{e}$ for each Ising vector $e$ and showed that each axis of the monstrous Griess algebra in 
\cite{C} is essentially a half of an Ising vector $e$ and $\tau_{e}$ is a $2A$-involution of the Monster simple group $\mathbb{M}$ constructed by Griess
\cite{G}. It has been proved by \cite{C} that the conjugacy class of the product of two $2A$-involutions $\tau\tau'$ is one of the nine classes $1A, 2A, 3A, 4A, 5A, 6A, 4B, 2B$ and $3C$ in  $\mathbb{M}$ and the inner product of the axis $e_{\tau}$, $e_{\tau'}$ is uniquely determined by the conjugacy
class. The above result is listed in terms of Ising vector and $\tau$-involution as follows: 
\begin{center}
\begin{tabular}{|c|c|c|c|c|c|c|c|c|c|}
\hline
$\left\langle \tau_{e}\tau_{f}\right\rangle ^{\mathbb{M}}$  & $1A$  & $2A$  & $3A$  & $4A$  & $5A$  & $6A$  & $4B$  & $2B$  & $3C$\tabularnewline
\hline
$\left\langle e,f\right\rangle $  & $1/4$  & $1/2^{5}$  & $13/2^{10}$  & $1/2^{7}$  & $3/2^{9}$  & $5/2^{10}$  & $1/2^{8}$  & $0$  & $1/2^{8}$\tabularnewline
\hline
\end{tabular}
\par\end{center}
Not only for the Ising vectors in the monstrous Griess algebra, In \cite{S}, Sakuma showed the inner product of two Ising vectors in any Moonshine type 
VOA is exactly one of the nine cases above. 

What really interests me is the subVOA structure generated by two Ising vectors in any moonshine type VOA. In \cite{LYY}, \cite{LYY1}, they studyed certain coset subalgebras $U$ of the lattice vertex operator algebra $V_{\sqrt{2}E_{8}}$, associated to sublattices $L$ of the $E_{8}$ lattice obtained by removing one node from the extended $E_{8}$ diagram at each time. In each of nine cases, $U$ contains two Ising vectors such that $U$ is generated by these two Ising vectors and their inner product corresponds to the values given in the table above. For the VOAs they constructed, the uniqueness of VOA structure for the $6A$ case was proved in \cite{DJY}, other cases except for $4A$ were discussed in \cite{DZ}. Now the question is: whether or not 
any subVOA generated by two Ising vectors in any moonshine type VOA is isomorphic to one of the nine case in \cite{LYY}, \cite{LYY1}. For the $3A$ case, i.e. $\langle e, f\rangle=\frac{13}{2^{10}}$, the answer is yes and given in \cite{M2} and \cite{SY}. In \cite{M}, Miyamoto also discussed the case $\langle e, f\rangle=\frac{1}{2^{8}}$ but didn't give a complete answer. Now in this paper, first I will finish the case $\langle e, f\rangle=\frac{1}{2^{8}}$ and then give a positive answer for the case $\langle e, f\rangle=\frac{3}{2^{9}}$.

The paper is organized as follows. In Section 2, we review some basic
notions and some results needed later. Most of notations follow the paper \cite{S}. In Section 3, we study the VOAs generated by two Ising vectors whose inner product is $\frac{1}{2^8}$ or $\frac{3}{2^9}$ and show that they are isomorphic to the cases constructed in \cite{LYY}, \cite{LYY1}.
Then from \cite{DZ}, we can immediately get that they both 
have unique VOA structures.
\section{Preliminary}
In this paper, we will use the same notation as in \cite{S}.
Let $V=(V,Y,\textbf{1},\omega)$ be a VOA over the real number field $\mathbb{R}$ which satisfies the following conditions:
\[
V=\sum_{n=0}^{\infty}V_{n}, V_{0}=\mathbb{R}\textbf{1}, V_{1}=0.
\]
Then by \cite{L}, there is a unique symmetric invariant bilinear form $\langle \ , \ \rangle$ on V such that $\langle \textbf{1},\textbf{1} \rangle=1$. We also assume that  $\langle \ , \ \rangle$ is positive definite. 
\begin{remark}\label{simpleremark}
Here the reason why I assume $\langle \ , \ \rangle$ is positive definite on $V$ is to guarantee $\langle \ , \ \rangle$ is nondegenerate on any subVOA of $V$.
Then by \cite{L}, any subVOA of $V$ is simple.
\end{remark}
\subsection{Ising vector and $\tau$-involution}

\begin{definition} A vector $e\in V_{2}$ is called a \emph{conformal
vector with central charge $c_{e}$ }if it satisfies\emph{ $e_{1}e=2e$
}and $e_{3}e=\frac{c_{e}}{2}\mathbf{1}$. Then the operators $L_{n}^{e}:=e_{n+1},\ n\in\mathbb{Z}$,
satisfy the Virasoro commutation relation
\[
\left[L_{m}^{e},\ L_{n}^{e}\right]=\left(m-n\right)L_{m+n}^{e}+\delta_{m+n,\ 0}\frac{m^{3}-m}{12}c_{e}
\]
for $m,\ n\in\mathbb{Z}.$ A conformal vector $e\in V_{2}$ with central
charge $\frac{1}{2}$ is called an \emph{Ising vector }if $e$ generates the
Virasoro vertex operator algebra $L(\frac{1}{2}, 0)$.

\end{definition}

Let $e$ be an Ising vector. By the definition, $e$ generates the Virasoro vertex operator algebra $L(\frac{1}{2}, 0)$. By \cite{DMZ}, $L(\frac{1}{2}, 0)$ is rational and has three irreducible modules $L(\frac{1}{2}, 0)$, $L(\frac{1}{2}, \frac{1}{2})$, $L(\frac{1}{2}, \frac{1}{16})$. Then as an  $L(\frac{1}{2}, 0)$-module, we have the following decomposition:
\[
V=V_{e}(0)\oplus V_{e}(\frac{1}{2})\oplus V_{e}(\frac{1}{16}).
\]
Here $V_{e}(h)$, $h=0, \frac{1}{2}, \frac{1}{16}$ is the sum of submodules isomorphic to $L(\frac{1}{2}, h)$.
Define $\tau_{e}$ in the following way:
\[
\tau_{e}(v)=
\begin{cases}
1& v\in V_{e}(0)\oplus V_{e}(\frac{1}{2}), \\
-1& v \in V_{e}(\frac{1}{16}).
\end{cases}
\]
Then by the fusion rules for $L(\frac{1}{2}, 0)$, $\tau_{e}$ is an automorphism of $V$ of order 2. We call it $\tau$-involution.

\subsection{Griess algebra $V_{2}$}
For any two elements $x, y \in V_{2}$, if we define the product $xy:=x_{1}y$, then $V_{2}$ becomes a commutative nonassociative algebra. Besides, $V_{2}$ has a bilinear form which is the restriction of $\langle \ , \ \rangle$ on $V_{2}$. Furthermore, for any $x, y, z \in V_{2}$, we have $\langle x, y\rangle\textbf{1}=x_{3}y$, $\langle xy, z \rangle=\langle y, xz\rangle$.
\begin{remark}\label{remarkfixed}
From the relation $\langle xy, z \rangle=\langle y, xz\rangle$, we can easily deduce that for any Ising vector $e$, $\langle \tau_{e}(x), \tau_{e}(y)\rangle=\langle x, y \rangle$ for any $x, y \in V_{2}$.
\end{remark}
The following lemma will be needed later:
\begin{lemma}\cite{M2}\label{miyamotodecomposition}
Let $e$ be an Ising vector. Then $V_{2}$ decomposes into
\[
V_{2}=\mathbb{R}e\oplus E^{e}(0)\oplus E^{e}(\frac{1}{2})\oplus E^{e}(\frac{1}{16}),
\]
where $E^{e}(h)$ denotes the eigenspace of $e_{1}$ with eigenvalue $h$.
\end{lemma}

\subsection{Subalgebra of $V_{2}$ generated by two Ising vectors $e$ and $f$}
Let $e, f$ be two Ising vectors. Now we consider the subalgebra of $V_{2}$ generated by $e$ and $f$. For any automorphsim $\sigma$ of V, we use $e^{\sigma}$ to denote the action of $\sigma$ on $e$. For any two elements $x, y \in V_{2}$, define $\alpha(x, y):=xy-\frac{1}{16}(x+y)$. In $\cite{S}$, Sakuma gives the following result:
\begin{lemma}\cite{S}\label{sakumaspan}
Let $X$ be the subalgebra generated by $e$ and $f$. Then $X$ is spanned by:
\[
S:=\{e, e^{\tau_{f}}, e^{\tau_{f}\tau_{e}}, f, f^{\tau_{e}}, f^{\tau_{e}\tau_{f}}, \alpha(e, f), \alpha(e, e^{\tau_{f}})\}.
\]

\end{lemma}
\begin{remark}\label{remarkfixed2}
$\alpha(e, f)=\alpha(f, e)$ and it is fixed by both $\tau_{e}$ and $\tau_{f}$.
\end{remark}
Let $T$ be the subgroup of $\mathrm{Aut}(V)$ generated by $\tau_{e}$ and $\tau_{f}$, $x^{T}$ denote the orbit of $x \in V_{2}$ under the action of $T$. Let $\rho=\tau_{e}\tau_{f}$. Then:
\begin{lemma}\cite{S}\label{Sakumalemma}\\
(1) $|e^{T}|=|f^{T}|$. In particular, $e=e^{\rho^n}$ if and only if $f=f^{\rho^n}$.\\
(2) $e^{T}=f^{T}$ if and only if $e^{T}$ is odd and $f=e^{\rho^{\frac{n+1}{2}}}$, where $n=|e^{T}|$.\\
(3) $(\tau_{e}\tau_{f})^{|e^{T}\bigcup f^{T}|}=1$ as an automorphism of $V$. 
\end{lemma}

\begin{theorem}\cite{S}
Let $N=|e^{T}\bigcup f^{T}|$.\\
(1) If $N=2$, then $\langle e, f\rangle=0$ or $\frac{1}{2^5}$.\\
(2) If $N=3$, then $\langle e, f\rangle=\frac{13}{2^{10}}$ or $\frac{1}{2^8}$.\\
(3) If $N=4$, then $(\langle e, f\rangle, \langle e, e^{\tau_{f}}\rangle)=(\frac{1}{2^7}, 0)$ or $(\frac{1}{2^8}, \frac{1}{2^5})$.\\
(4) If $N=5$, then $\langle e, f\rangle=\langle e, e^{\tau_{f}}\rangle=\frac{3}{2^9}$.\\
(5) If $N=6$, then $\langle e, f\rangle=\frac{5}{2^{10}}$, $\langle e, e^{\tau_{f}}\rangle=\frac{13}{2^{10}}$ and $\langle e^{\tau_{f}}, f^{\tau_{e}}\rangle=\frac{1}{2^5}$.
\end{theorem}

\subsection{Some calculations in $V_{2}$}\label{calculation}
Let $e, f$ be two Ising vectors. Then $\langle e, e \rangle=\langle f, f \rangle=\frac{1}{4}$ by definition. Let $a$ be any Ising vector, Denote $\langle e, f \rangle:=\frac{\lambda_{1}}{4}$, 
$\langle e, e^{\tau_{f}}\rangle:=\frac{\lambda_{2}}{4}$. Here I list some results computed in $\cite{S}$:
\begin{align}
 & a\cdot\alpha(a, x)=\frac{7}{16}\alpha(a, x)+(12\langle a, x\rangle-\frac{25}{2^8})a+\frac{7}{2^9}(x+x^{\tau_{a}}), x\in V_{2}\\ 
& \langle a, \alpha(a, x)\rangle=\frac{31}{16}\langle a, x\rangle-\frac{1}{2^6}, x\in V_{2}
\end{align}

\begin{align}
\alpha(e, f)\cdot\alpha(e, f) 
& =\frac{7}{3}(4\lambda_{1}^2-\frac{1}{2^4}\lambda_{1}-\frac{1}{2^{12}}+\frac{1}{6}\lambda_{2})e+
\frac{7^2}{3\cdot 2^5}(\lambda_{1}-\frac{5}{2^8})(f+f^{\tau_{e}})\nonumber\\
& +\frac{7^2}{3\cdot 2^{13}}(e^{\tau_{f}}+e^{\tau_{f}\tau_{e}})-\frac{1}{3}(5\lambda_{1}+\frac{13}{2^7})\alpha(e, f)\nonumber\\
& -\frac{7}{3\cdot 2^7}\alpha(e, e^{\tau_{f}})+ \frac{7}{2^9}\alpha(f, f^{\tau_{e}}),
\end{align}

\begin{align}
\frac{1}{7}(2^{11}\lambda_{1}^2
& -9\cdot 2^{4}\lambda_{1}+\frac{33}{2^4}+2^3\lambda_{2})(e-f)+ (2^4\lambda_{1}-\frac{3}{8})(f^{\tau_{e}}-e^{\tau_{f}})\nonumber\\
& +\frac{1}{2^4}(e^{\tau_{f}\tau_{e}}-f^{\tau_{e}\tau_{f}})-(\alpha(e, e^{\tau_{f}})-\alpha(f, f^{\tau_{e}}))=0. 
\end{align}
We can compute $\langle \alpha(e, f), \alpha(e, f)\rangle$ in the following way:
\begin{align*}
\langle \alpha(e, f), \alpha(e, f)\rangle
& =\langle ef-\frac{1}{16}(e+f), \alpha(e, f)\rangle \\
& =\langle ef, \alpha(e, f)\rangle-\frac{1}{16}\langle e, \alpha(e, f)\rangle-\frac{1}{16}\langle f, \alpha(e, f)\rangle\\
& =\langle f, e\alpha(e, f)\rangle-\frac{1}{16}\langle e, \alpha(e, f)\rangle-\frac{1}{16}\langle f, \alpha(e, f)\rangle\\
& =\langle f, \frac{7}{16}\alpha(e,f)+(12\langle e, f\rangle-\frac{25}{2^8})e+\frac{7}{2^9}(f+f^{\tau_{e}})\rangle-\frac{1}{16}\langle e, \alpha(e, f)\rangle-\frac{1}{16}\langle f, \alpha(e, f)\rangle\\
& =\frac{6}{16}\langle f, \alpha(e, f)\rangle+(12\langle e, f\rangle-\frac{25}{2^8})\langle f, e\rangle+\frac{7}{2^9}\langle f, f\rangle
+\frac{7}{2^9}\langle f, f^{\tau_{e}}\rangle-\frac{1}{16}\langle e, \alpha(e, f)\rangle.
\end{align*}
For $\langle f, \alpha(e, e^{\tau_{f}})\rangle$, by remark \ref{remarkfixed} and \ref{remarkfixed2}, we have
\begin{align*}
\langle f, \alpha(e, e^{\tau_{f}})\rangle
=\langle f^{\tau_{e}\tau_{f}}, \alpha(e, f)\rangle,
\end{align*}
If $f^{\tau_{e}\tau_{f}}=e^{\tau_{f}\tau_{e}}$, then the above equation becomes
\begin{align*}
\langle f, \alpha(e, e^{\tau_{f}})\rangle
=\langle f^{\tau_{e}\tau_{f}}, \alpha(e, f)\rangle=\langle e^{\tau_{f}\tau_{e}}, \alpha(e, f)\rangle=\langle e, \alpha(e, f)\rangle.
\end{align*}

\section{SubVOA generated by two Ising vectors $e$ and $f$}
We denote the subVOA generated by $e$ and $f$ as $\mathrm{VA}(e,f)$.

\subsection{The case $\langle e, f\rangle=\frac{1}{2^8}$}
By Lemma \ref{miyamotodecomposition}, we can write
\[
f=\frac{1}{2^6}e+a+b+c,
\]
where $a\in E^{e}(0)$, $b\in E^{e}(\frac{1}{2})$, $c\in E^{e}(\frac{1}{16})$. The following three lemmas were given in $\cite{M2}$:
\begin{lemma}\label{C3eac}
$b=0$, $\mathrm{VA}(e,f)_{2}$ has a basis $\{e, a, c\}$.
\end{lemma}

\begin{remark}
If we use the language in lemma $\ref{sakumaspan}$, then we have $\mathrm{VA}(e,f)_{2}$ has a basis $\{e, f, e^{\tau_{f}}\}$.
\end{remark}

\begin{lemma}\label{C3Virasoro}
Let $\widetilde{\omega}=\frac{64}{33}a$. Then $e$ and $\widetilde{\omega}$ are two orthogonal Ising vectors with central charge $\frac{1}{2}$, $\frac{21}{22}$ 
respectively and $e+\widetilde{\omega}$ is a Virasoro element of $\mathrm{VA}(e,f)$.
\end{lemma}

For simplicity, we will write $[h_{1}, h_{2}]$ for $L(\frac{1}{2}, h_{1})\otimes L(\frac{21}{22}, h_{2})$
\begin{lemma}\label{C3multiplicity}
As a $L(\frac{1}{2}, 0)\otimes L(\frac{21}{22}, 0)$-module,
\[
\mathrm{VA}(e,f)=n_{1}[0, 0]+n_{2}[0, 8]+n_{3}[\frac{1}{2}, \frac{7}{2}]+n_{4}[\frac{1}{2}, \frac{45}{2}]+n_{5}[\frac{1}{16}, \frac{31}{16}]
+n_{6}[\frac{1}{16}, \frac{175}{16}].
\]
where $n_{1}=1, n_{i}, i=2,3,4,5,6$ are nonnegative integers.
\end{lemma}

Now we can state our result:
\begin{theorem}\label{C3structure}
As a $L(\frac{1}{2}, 0)\otimes L(\frac{21}{22}, 0)$-module,
\begin{eqnarray*}
\mathrm{VA}(e,f)=[0, 0]+[0, 8]+[\frac{1}{2}, \frac{7}{2}]+[\frac{1}{2}, \frac{45}{2}]+[\frac{1}{16}, \frac{31}{16}]
+[\frac{1}{16}, \frac{175}{16}],
\end{eqnarray*}
i.e. all $n_{i}, i=1,2,3,4,5,6$ in lemma $\ref{C3multiplicity}$ equal to $1$.
\end{theorem} 
 
\begin{proof}
\begin{claim}
$n_{5}=1$.
\end{claim}
From lemma $\ref{C3eac}$ and lemma $\ref{C3Virasoro}$, we have $e_{1}c=\frac{1}{16}c$, $\widetilde{\omega}_{1}c=\frac{31}{16}c$. Since $\widetilde{\omega}_{2}c \in V_{1}$ which is $0$, $\widetilde{\omega}_{3}c=\langle \widetilde{\omega}, c\rangle\textbf{1}=0$, $\widetilde{\omega}_{n}c=0$ for $n>3$, we have $c$ is a highest weight vector of $L(\frac{21}{22}, 0)$ with highest weight $\frac{31}{16}$. So $[\frac{1}{16}, \frac{31}{16}]$ appears in $\mathrm{VA}(e,f)$. Then the Claim follows from lemma $\ref{C3eac}$.

\begin{claim}
As a $L(\frac{1}{2}, 0)\otimes L(\frac{21}{22}, 0)$-module, $[\frac{1}{2}, \frac{7}{2}] \subseteq \mathrm{VA}(e,f)$.
\end{claim}
First, from $\cite{M2}$, we have the following results:
\[
ac=\frac{1023}{2^{10}}c,\ \ \ \ \ cc=\frac{63}{2^{11}}e+\frac{31}{32}a,\ \ \ \ \ \ \langle c, c\rangle=\frac{63}{2^9}.
\]
Then we compute $\langle c_{-1}c, c_{-1}c \rangle$:
\begin{align}\label{c-1c}
\langle c_{-1}c, c_{-1}c \rangle & =\langle c, c_{3}c_{-1}c \rangle=\langle c, [c_{3}, c_{-1}]c \rangle+\langle c, c_{-1}c_{3}c \rangle \nonumber\\
 & = \langle c, (c_{0}c)_{2}c+3(c_{1}c)_{1}c+3(c_{2}c)_{0}c+(c_{3}c)_{-1}c\rangle+\langle c, c\rangle^{2} \nonumber\\
& =\langle c, ([c_{1}, c_{1}]-(c_{1}c)_{1})c+3(c_{1}c)_{1}c+(c_{3}c)_{-1}c\rangle+\langle c, c\rangle^{2} \nonumber\\
& =\langle c, 2(c_{1}c)_{1}c\rangle+2\langle c, c\rangle^{2} \nonumber\\
& =\frac{63}{2^{10}}\langle c, ec\rangle+ \frac{31}{16}\langle c, ac\rangle+ 2\langle c, c\rangle^{2} \nonumber\\
& =\frac{63}{2^{14}}\langle c, c\rangle+ \frac{31}{16}\cdot \frac{1023}{2^{10}}\langle c, c\rangle+ 2\langle c, c\rangle^{2} \nonumber\\
& =\frac{1119\cdot 63}{2^{18}} \neq 0,
\end{align}
So $c_{-1}c \neq 0$.

From the fusion rules for $L(\frac{1}{2}, 0)$ and $L(\frac{21}{22}, 0)$, we have 
\begin{eqnarray}\label{1/16cross1/16}
c_{-1}c \in[\frac{1}{16}, \frac{31}{16}]. [\frac{1}{16}, \frac{31}{16}] \subseteq [0, 0]+ [0, 8]+ [\frac{1}{2}, \frac{7}{2}].
\end{eqnarray}
Assume $[\frac{1}{2}, \frac{7}{2}]$ does not appear in $\mathrm{VA}(e,f)$.
Since $\mathrm{wt}c_{-1}c=4$, Then by ($\ref{1/16cross1/16}$),
\begin{eqnarray*}
c_{-1}c \in [0, 0].
\end{eqnarray*}
Notice that the weight $4$ subspace of $[0, 0]$ has a basis $\{e_{-3}\textbf{1}, \widetilde{\omega}_{-3}\textbf{1}, e_{-1}e,
 \widetilde{\omega}_{-1}\widetilde{\omega}, e_{-1}\widetilde{\omega} \}$. Set 
\begin{align*}
& x_{1}=e_{-3}\textbf{1}, \ \ \ \ x_{2}=\widetilde{\omega}_{-3}\textbf{1},\ \ \ \  x_{3}=e_{-1}e,\\
& x_{4}=\widetilde{\omega}_{-1}\widetilde{\omega},\ \ \ \ x_{5}=e_{-1}\widetilde{\omega}, \ \ \ \ x_{6}=c_{-1}c,
\end{align*}
By a similar way as ($\ref{c-1c}$), we can compute all $\langle x_{i}, x_{j}\rangle, 1\leqslant i, j \leqslant 6$. Then we obtain the matrix
$(\langle x_{i}, x_{j}\rangle)_{1\leqslant i, j \leqslant 6}$:
\begin{equation}
\left(
\begin{array}{cccccc}
\frac{5}{2} & 0 & \frac{3}{2} & 0 & 0 & \frac{3\cdot 63}{2^{13}}\\[0.2cm]
0 & \frac{5}{2} & 0 & \frac{3}{2} & 0 & \frac{63\cdot 31\cdot 3}{2^{13}}\\[0.2cm]
\frac{3}{2} & 0 & \frac{17}{8} & 0 & 0 & \frac{63\cdot 33}{2^{17}}\\[0.2cm]
0 & \frac{3}{2} & 0 & \frac{17}{8} & 0 & \frac{31\cdot 33 \cdot 63}{2^{17}}\\[0.2cm]
0 & 0 & 0 & 0 & \frac{1}{16} & \frac{31\cdot 33}{2^{17}}\\[0.2cm]
\frac{3\cdot 63}{2^{13}} & \frac{63\cdot 31\cdot 3}{2^{13}} & \frac{63\cdot 33}{2^{17}} & \frac{31\cdot 33 \cdot 63}{2^{17}} &
\frac{31\cdot 33}{2^{17}} & \frac{1119\cdot 63}{2^{18}}\\[0.2cm]
\end{array}
\right).
\end{equation}
The $\mathrm{det}(\langle x_{i}, x_{j}\rangle)_{1\leqslant i, j \leqslant 6}=0.0093 \neq 0$, which implies $\{x_{i}, i=1, ..., 6\}$ is linearly independent,
contradicting with the fact that $\{x_{i}, i=1, ..., 5\}$ is a basis of wight $4$ subspace of $[0, 0]$. Hence $[\frac{1}{2}, \frac{7}{2}] \subseteq  \mathrm{VA}(e,f)$.
\begin{claim}
$n_{2}=n_{3}$, $n_{4}=1$.
\end{claim}
Let
\begin{align*}
U=[0, 0]+n_{2}[0, 8]+n_{3}[\frac{1}{2}, \frac{7}{2}]+n_{4}[\frac{1}{2}, \frac{45}{2}],
\end{align*}
By the fusion rules for $L(\frac{1}{2}, 0)$, we have $U$ is a subVOA of $\mathrm{VA}(e,f)$. By Remark $\ref{simpleremark}$, $U$ is simple. By $\cite{ABD}$
and $\cite{HKL}$, $U$ is rational and $C_{2}$-cofinite. Define $\sigma$ as follows:
\[
\sigma(v)=
\begin{cases}
1& v\in [0, 0]+n_{2}[0, 8], \\
-1& v \in n_{3}[\frac{1}{2}, \frac{7}{2}]+n_{4}[\frac{1}{2}, \frac{45}{2}].
\end{cases}
\]
Since $[\frac{1}{2}, \frac{7}{2}] \subseteq  \mathrm{VA}(e,f)$, then by the fusion rules for $L(\frac{1}{2}, 0)$ again we have $\sigma$ is an automorphism of $U$ with order $2$. From Quantum Galois theory $\cite{DM2}$ and $\cite{ADJR}$ as well as quantum dimension \cite{DJX},
 we have

\begin{eqnarray}
q\dim_{[0, 0]+n_{2}[0, 8]}(n_{3}[\frac{1}{2}, \frac{7}{2}]+n_{4}[\frac{1}{2}, \frac{45}{2}])=1,
\end{eqnarray}
equivalently,
\begin{eqnarray}\label{qdimc31}
q\dim_{[0, 0]}([0, 0]+n_{2}[0, 8])=q\dim_{[0, 0]}(n_{3}[\frac{1}{2}, \frac{7}{2}]+n_{4}[\frac{1}{2}, \frac{45}{2}]).
\end{eqnarray}
By direct computation, we have
\begin{align*}
& q\dim_{[0, 0]}([0, 8])=q\dim_{[0, 0]}([\frac{1}{2}, \frac{7}{2}])=2+\sqrt{3},\\
& q\dim_{[0, 0]}([0, 0])=q\dim_{[0, 0]}([\frac{1}{2}, \frac{45}{2}])=1.
\end{align*}
Then ($\ref{qdimc31}$) implies 
\begin{align*}
1+(2+\sqrt{3})n_{2}=(2+\sqrt{3})n_{3}+n_{4},
\end{align*}
Hence $n_{2}=n_{3}$, $n_{4}=1$.

\begin{claim}
$n_{6}=1$.
\end{claim}
Since $q\dim_{[0, 0]}([\frac{1}{2}, \frac{45}{2}])=1$, by \cite{DJX} we have $[\frac{1}{2}, \frac{45}{2}]$ is a simple current module of 
$[0, 0]$. By the fusion rules for $[0, 0]$-modules  (see \cite{W} the formula of the fusion rules for minimal models of Virasoro algebra) we have
\begin{align*}
& [\frac{1}{2}, \frac{45}{2}]\boxtimes [\frac{1}{16}, \frac{31}{16}]=[\frac{1}{16}, \frac{175}{16}],\\
& [\frac{1}{2}, \frac{45}{2}]\boxtimes [\frac{1}{16}, \frac{175}{16}]=[\frac{1}{16}, \frac{31}{16}].
\end{align*}
Since $n_{4}=n_{5}=1$, a similar proof as proposition 5.1(1) in $\cite{DMZ}$ shows that $n_{6}$ must be $1$.

\begin{claim}
$n_{2}=n_{3}=1$.
\end{claim}
Let
\begin{align*}
& U^{1}=[0, 0]+n_{2}[0, 8]+n_{3}[\frac{1}{2}, \frac{7}{2}]+[\frac{1}{2}, \frac{45}{2}],\\
& U^{2}=[\frac{1}{16}, \frac{31}{16}]+[\frac{1}{16}, \frac{175}{16}],
\end{align*}
Then by the fusion rules for $L(\frac{1}{2}, 0)$ and Claims 1-4, we have 
\begin{align*}
&\mathrm{VA}(e,f)=U^{1}+U^{2},\\
& U^{1}.U^{2}\subseteq U^{2}, \ \ \ U^{2}.U^{1}\subseteq U^{2},\ \ \ U^{2}.U^{2}\subseteq U^{1}.
\end{align*}
Define $\sigma$ as follows:
\[
\sigma(v)=
\begin{cases}
1& v\in U^{1}, \\
-1& v \in U^{2}.
\end{cases}
\]
Same reason as Claim 3, we have $\sigma$ is an automorphism of $\mathrm{VA}(e,f)$ with order $2$, and 
\begin{eqnarray}\label{qdimc32}
q\dim_{[0, 0]}(U^{1})=q\dim_{[0, 0]}(U^{2}).
\end{eqnarray}
By direct computation, we have 
\begin{align*}
q\dim_{[0, 0]}([\frac{1}{16}, \frac{31}{16}])=q\dim_{[0, 0]}([\frac{1}{16}, \frac{175}{16}])=3+\sqrt{3}.
\end{align*}
Then ($\ref{qdimc32}$) implies 
\begin{align*}
1+(2+\sqrt{3})n_{2}+(2+\sqrt{3})n_{3}+1=3+\sqrt{3}+3+\sqrt{3},
\end{align*}
From Claim 3, we have $n_{2}=n_{3}$. Combining it with the equation above we can get $n_{2}=n_{3}=1$.
\end{proof}

\begin{theorem}
When $\langle e, f\rangle=\frac{1}{2^8}$, $\mathrm{VA}(e,f)$ has a unique VOA structure.
\end{theorem}
\begin{proof}
The theorem follows directly from theorem 4.6 in $\cite{DZ}$ and theorem $\ref{C3structure}$ above.
\end{proof}
\subsection{The case $\langle e, f\rangle=\frac{3}{2^9}$}
\begin{remark}\cite{S}
In this case, $\langle e, e^{\tau_{f}}\rangle$ is also equal to $\frac{3}{2^9}$.
\end{remark}

First we will focus on the weight 2 subspace of $\mathrm{VA}(e,f)$.
\begin{lemma}
The dimension of $\mathrm{VA}(e,f)_{2}$ is six and it has a basis $\{e, e^{\tau_{f}}, e^{\tau_{f}\tau_{e}}, f, f^{\tau_{e}}, \alpha(e, f)\}$.
\end{lemma}
\begin{proof}
From lemma $\ref{Sakumalemma}$, we have $e^{\tau_{f}\tau_{e}}=f^{\tau_{e}\tau_{f}}$. Then by lemma $\ref{sakumaspan}$, we can get 
\begin{align*}
\mathrm{VA}(e,f)_{2}=\mathrm{span}\{e, e^{\tau_{f}}, e^{\tau_{f}\tau_{e}}, f, f^{\tau_{e}}, \alpha(e, f), \alpha(e, e^{\tau_{f}})\}.
\end{align*}
Let
\begin{align*}
& v_{1}=e, \ \ \ \ v_{2}=e^{\tau_{f}},\ \ \ \ v_{3}=e^{\tau_{f}\tau_{e}}, \ \ \ \ v_{4}=f,\\
& v_{5}=f^{\tau_{e}},\ \ \ \ v_{6}=\alpha(e, f), \ \ \ \ v_{7}=\alpha(e, e^{\tau_{f}}).
\end{align*}
Then by using results in section $\ref{calculation}$, we can find all $\langle v_{i}, v_{j}\rangle, 1\leq i, j \leq 7$. Let
 $A=(\langle v_{i}, v_{j}\rangle)_{1\leq i, j \leq 7}$. Then:
\begin{equation}
A=\left(
\begin{array}{ccccccc}
\frac{1}{4} & \frac{3}{2^9} & \frac{3}{2^9} & \frac{3}{2^9} & \frac{3}{2^9} & \frac{-35}{2^{13}} & \frac{-35}{2^{13}}\\[0.2cm]
\frac{3}{2^9} & \frac{1}{4} & \frac{3}{2^9} & \frac{3}{2^9} & \frac{3}{2^9} & \frac{-35}{2^{13}} & \frac{-35}{2^{13}}\\[0.2cm]
\frac{3}{2^9} & \frac{3}{2^9} & \frac{1}{4} & \frac{3}{2^9} & \frac{3}{2^9} & \frac{-35}{2^{13}} & \frac{-35}{2^{13}}\\[0.2cm]
\frac{3}{2^9} & \frac{3}{2^9} & \frac{3}{2^9} & \frac{1}{4} & \frac{3}{2^9} & \frac{-35}{2^{13}} & \frac{-35}{2^{13}}\\[0.2cm]
\frac{3}{2^9} & \frac{3}{2^9} & \frac{3}{2^9} & \frac{3}{2^9} & \frac{1}{4} & \frac{-35}{2^{13}} & \frac{-35}{2^{13}}\\[0.2cm]
\frac{-35}{2^{13}} & \frac{-35}{2^{13}} & \frac{-35}{2^{13}} & \frac{-35}{2^{13}} & \frac{-35}{2^{13}} & \frac{525}{2^{18}} & \frac{-175}{2^{17}}\\[0.2cm]
\frac{-35}{2^{13}} & \frac{-35}{2^{13}} & \frac{-35}{2^{13}} & \frac{-35}{2^{13}} & \frac{-35}{2^{13}} & \frac{-175}{2^{17}} & \frac{525}{2^{18}}\\[0.2cm]
\end{array}
\right).
\end{equation}
Let $AX^{T}=0$, where $X=(x_{1}, x_{2}, ... ,x_{7})$. By solving this linear systems, we can find the following relation:
\begin{align*}
x_{1}=x_{2}=x_{3}=x_{4}=x_{5}=\frac{1}{2^5}x_{6}=\frac{1}{2^5}x_{7}.
\end{align*}
So $\mathrm{rank}(A)=6$ and 
\begin{align}\label{A5basisrelation}
e+e^{\tau_{f}}+e^{\tau_{f}\tau_{e}}+f+f^{\tau_{e}}+2^5 \alpha(e, f)+2^5\alpha(e, e^{\tau_{f}})=0.
\end{align}
This completes the proof.
\end{proof}
Our next goal is to find mutually orthogonal conformal vectors.
\begin{lemma}\label{A5lemmabasis}
Let 
\begin{align*}
& u=-\frac{1}{2^3}e+f+f^{\tau_{e}}-\frac{2^5}{7}\alpha(e, f),\\
& v=\frac{1}{2^5}e+e^{\tau_{f}}+e^{\tau_{f}\tau_{e}}+\frac{2^5}{7}\alpha(e, f).
\end{align*}
Then $\{u, v\}$ is a basis of $E^{e}(0)$.
\end{lemma}
\begin{proof}
First we need to do some calculations:
\begin{align*}
& ee=2e,\\
& ee^{\tau_{f}}=\alpha(e, e^{\tau_{f}})+\frac{1}{16}(e+e^{\tau_{f}})
=\frac{1}{2^5}e+\frac{1}{2^5}e^{\tau_{f}}-\frac{1}{2^5}e^{\tau_{f}\tau_{e}}-\frac{1}{2^5}f-\frac{1}{2^5}f^{\tau_{e}}-
\alpha(e, f),\\
& ee^{\tau_{f}\tau_{e}}=\tau_{e}(ee^{\tau_{f}})
=\frac{1}{2^5}e-\frac{1}{2^5}e^{\tau_{f}}+\frac{1}{2^5}e^{\tau_{f}\tau_{e}}-\frac{1}{2^5}f-\frac{1}{2^5}f^{\tau_{e}}-
\alpha(e, f),\\
& ef=\frac{1}{2^4}e+\frac{1}{2^4}f+\alpha(e, f),\\
& ef^{\tau_{e}}=\tau_{e}(ef)=\frac{1}{2^4}e+\frac{1}{2^4}f^{\tau_{e}}+\alpha(e, f),\\
& e\alpha(e, f)=-\frac{7}{2^8}e+\frac{7}{2^9}f+\frac{7}{2^9}f^{\tau_{e}}+\frac{7}{2^4}\alpha(e, f).
\end{align*}
For the last equation above, see section $\ref{calculation}$. Now consider the following equation:
\begin{equation}\label{A5e0basis}
e\cdot(x_{1}e+x_{2}e^{\tau_{f}}+x_{3}e^{\tau_{f}\tau_{e}}+x_{4}f+x_{5}f^{\tau_{e}}+x_{6}\alpha(e, f))=0.
\end{equation}
($\ref{A5e0basis}$) gives us $AX^{T}=0$, where
\begin{equation}\nonumber
A=\left(
\begin{array}{cccccc}
2 & \frac{1}{2^5} & \frac{1}{2^5} & \frac{1}{2^4} & \frac{1}{2^4} & -\frac{7}{2^{8}}\\[0.2cm]
0 & \frac{1}{2^5} & -\frac{1}{2^5} & 0 & 0 & 0 \\[0.2cm]
0 & -\frac{1}{2^5} & \frac{1}{2^5} & 0 & 0 & 0\\[0.2cm]
0 & -\frac{1}{2^5} & -\frac{1}{2^5} & \frac{1}{2^4} & 0 & \frac{7}{2^{9}}\\[0.2cm]
0 & -\frac{1}{2^5} & -\frac{1}{2^5} & 0 & \frac{1}{2^4} & \frac{7}{2^{9}}\\[0.2cm]
0 & -1 & -1 & 1 & 1 & \frac{7}{2^{4}}\\[0.2cm]
\end{array}
\right),
\end{equation}
\[
X=(x_{1}, x_{2}, ... ,x_{6}).
\]
By solving it, we have 
\begin{align*}
& x_{1}=\frac{1}{2^5}x_{3}-\frac{1}{2^3}x_5,\ \ \ \ x_{2}=x_{3},\\
& x_{4}=x_{5},\ \ \ \ x_{6}=\frac{2^5}{7}x_{3}-\frac{2^5}{7}x_{5}.
\end{align*}
Then $x_{3}=0, x_{5}=1$ gives $u$, $x_{3}=1, x_{5}=0$ gives $v$.
\end{proof}
\begin{lemma}\label{A5conformalvectors}
Let 
\begin{align*}
 \widetilde{u}=\frac{112}{105k_{1}+140}(u+k_{1}v),\ \ \ 
 \widetilde{v}=\frac{112}{105k_{2}+140}(u+k_{2}v),
\end{align*}
where
\begin{align*}
k_{1}=\frac{124+56\sqrt{5}}{19}, \ \ \  k_{2}=\frac{124-56\sqrt{5}}{19}.
\end{align*}
Then $e, \widetilde{u}, \widetilde{v}$ are mutually orthogonal conformal vectors with central charge $\frac{1}{2}, \frac{25}{28}, \frac{25}{28}$
respectively.
\end{lemma}
\begin{proof}
Like what we did in section $\ref{calculation}$ and lemma $\ref{A5lemmabasis}$, we can find all $xy, \langle x, y \rangle$ for $x, y \in \{e, e^{\tau_{f}}, e^{\tau_{f}\tau_{e}}, f, f^{\tau_{e}}, \alpha(e, f)\}$. Then we can obtain the following:
\begin{align*}
&u^2=\frac{33}{14}u-\frac{1}{7}v,\ \ \ v^2=\frac{19}{112}u+\frac{229}{112}v,\ \ \ uv=-\frac{19}{112}u+\frac{1}{7}v,\\
&\langle u, u \rangle=\frac{1125}{7\cdot 2^8},\ \ \ \ \  \langle v, v \rangle=\frac{27250}{7\cdot 2^{13}},\ \ \ \ \ \ \ \ \
 \langle u, v \rangle=-\frac{125}{7\cdot 2^{10}}.
\end{align*}
Then it is a routine work to check that 
\begin{align*}
& \widetilde{u} \widetilde{u}=2 \widetilde{u},\ \ \ \widetilde{v}\widetilde{v}=2\widetilde{v},\ \ \ \widetilde{u}\widetilde{v}=0\\
&\langle \widetilde{u}, \widetilde{u} \rangle=\langle \widetilde{v}, \widetilde{v} \rangle=\frac{25}{56},\ \ \ 
 \langle \widetilde{u}, \widetilde{v} \rangle=0.
\end{align*}
Then we can complete our lemma by the equations above and the fact that $\widetilde{u}, \widetilde{v} \in E^{e}(0)$.
\end{proof}
For simplicity, we will write $[h_{1}, h_{2}, h_{3}]$ for $L(\frac{1}{2}, h_{1})\otimes L(\frac{25}{28}, h_{2})\otimes L(\frac{25}{28}, h_{2})$.
We have the following theorem:
\begin{theorem}\label{A5structure}
As a $L(\frac{1}{2}, h_{1})\otimes L(\frac{25}{28}, h_{2})\otimes L(\frac{25}{28}, h_{2})$-module,
\begin{align*}
\mathrm{VA}(e,f) 
& =[0, 0, 0]\oplus [0, \frac{15}{2}, \frac{15}{2}]\oplus [0, \frac{3}{4}, \frac{13}{4}]\oplus [0, \frac{13}{4}, \frac{3}{4}]\\
& \oplus[\frac{1}{2}, 0, \frac{15}{2}]\oplus [\frac{1}{2}, \frac{15}{2}, 0]\oplus [\frac{1}{2}, \frac{3}{4}, \frac{3}{4}]\oplus [\frac{1}{2}, \frac{13}{4}, \frac{13}{4}]\\
& \oplus[\frac{1}{16}, \frac{5}{32}, \frac{57}{32}]\oplus [\frac{1}{16}, \frac{57}{32}, \frac{5}{32}]\oplus [\frac{1}{16}, \frac{57}{32}, \frac{165}{32}]
\oplus [\frac{1}{16}, \frac{165}{32}, \frac{57}{32}].
\end{align*}
\end{theorem}
\begin{proof}
We write 
\begin{align*}
\mathrm{VA}(e,f) 
& =n_{1}[0, 0, 0]\oplus n_{2}[0, \frac{15}{2}, \frac{15}{2}]\oplus n_{3}[0, \frac{3}{4}, \frac{13}{4}]\oplus n_{4}[0, \frac{13}{4}, \frac{3}{4}]\\
& \oplus n_{5}[\frac{1}{2}, 0, \frac{15}{2}]\oplus n_{6}[\frac{1}{2}, \frac{15}{2}, 0]\oplus n_{7}[\frac{1}{2}, \frac{3}{4}, \frac{3}{4}]\oplus n_{8}[\frac{1}{2}, \frac{13}{4}, \frac{13}{4}]\\
& \oplus n_{9}[\frac{1}{16}, \frac{5}{32}, \frac{57}{32}]\oplus n_{10}[\frac{1}{16}, \frac{57}{32}, \frac{5}{32}]\oplus n_{11}[\frac{1}{16}, \frac{57}{32}, \frac{165}{32}] \oplus n_{12}[\frac{1}{16}, \frac{165}{32}, \frac{57}{32}].
\end{align*}
We will show that $n_{i}=1, i=1, ... ,12$.

\emph{Claim 1) $n_{1}=1$.}

The claim follows directly from lemma $\ref{A5conformalvectors}$.

\emph{Claim 2) $n_{7}=1$.}

Consider the following equation:
\begin{align}\label{A5e1/2basis}
& e\cdot(x_{1}e+x_{2}e^{\tau_{f}}+x_{3}e^{\tau_{f}\tau_{e}}+x_{4}f+x_{5}f^{\tau_{e}}+x_{6}\alpha(e, f)) \nonumber\\
& =\frac{1}{2}(x_{1}e+x_{2}e^{\tau_{f}}+x_{3}e^{\tau_{f}\tau_{e}}+x_{4}f+x_{5}f^{\tau_{e}}+x_{6}\alpha(e, f))
\end{align}
($\ref{A5e1/2basis}$) gives us $AX^{T}=0$, where
\begin{equation}\nonumber
A=\left(
\begin{array}{cccccc}
\frac{3}{2} & \frac{1}{2^5} & \frac{1}{2^5} & \frac{1}{2^4} & \frac{1}{2^4} & -\frac{7}{2^{8}}\\[0.2cm]
0 & -\frac{15}{2^5} & -\frac{1}{2^5} & 0 & 0 & 0 \\[0.2cm]
0 & -\frac{1}{2^5} & -\frac{15}{2^5} & 0 & 0 & 0\\[0.2cm]
0 & -\frac{1}{2^5} & -\frac{1}{2^5} & -\frac{7}{2^4} & 0 & \frac{7}{2^{9}}\\[0.2cm]
0 & -\frac{1}{2^5} & -\frac{1}{2^5} & 0 & -\frac{7}{2^4} & \frac{7}{2^{9}}\\[0.2cm]
0 & -1 & -1 & 1 & 1 & -\frac{1}{2^{4}}\\[0.2cm]
\end{array}
\right),
\end{equation}
\[
X=(x_{1}, x_{2}, ... ,x_{6}).
\]
By solving it, we have 
\begin{align*}
& x_{1}=\frac{1}{2^6}x_{6},\ \ \ \ x_{2}=x_{3}=0,\\
& x_{4}=x_{5}=\frac{1}{2^5}x_{6}.
\end{align*}
So $\dim E^{e}(\frac{1}{2})=1$ and it has a basis $\{a=e+2f+2f^{\tau{e}}+2^6\alpha(e, f)\}$. Besides, By direct computation, 
we have 
\begin{align*}
 u\cdot a=\frac{15}{2^4}a,\ \ \ \ v\cdot a=\frac{45}{2^6}a,
\end{align*}
Then
\begin{align*}
& \widetilde{u}\cdot a=\frac{112}{105k_{1}+140}(\frac{15}{2^4}+k_{1}\frac{45}{2^6})a=\frac{3}{4}a\\
& \widetilde{v}\cdot a=\frac{112}{105k_{2}+140}(\frac{15}{2^4}+k_{2}\frac{45}{2^6})a=\frac{3}{4}a.
\end{align*}
A similar argument as Claim 1 in theorem $\ref{C3structure}$ together with the fact $\dim E^{e}(\frac{1}{2})=1$ tells us
$n_{7}=1$

 \emph{Claim 3) $n_{9}=n_{10}=1$.}

Consider the following equation:
\begin{align}\label{A5e1/16basis}
& e\cdot(x_{1}e+x_{2}e^{\tau_{f}}+x_{3}e^{\tau_{f}\tau_{e}}+x_{4}f+x_{5}f^{\tau_{e}}+x_{6}\alpha(e, f)) \nonumber\\
& =\frac{1}{16}(x_{1}e+x_{2}e^{\tau_{f}}+x_{3}e^{\tau_{f}\tau_{e}}+x_{4}f+x_{5}f^{\tau_{e}}+x_{6}\alpha(e, f))
\end{align}
($\ref{A5e1/16basis}$) gives us $BX^{T}=0$, where
\begin{equation}\nonumber
B=\left(
\begin{array}{cccccc}
\frac{31}{2^4} & \frac{1}{2^5} & \frac{1}{2^5} & \frac{1}{2^4} & \frac{1}{2^4} & -\frac{7}{2^{8}}\\[0.2cm]
0 & -\frac{1}{2^5} & -\frac{1}{2^5} & 0 & 0 & 0 \\[0.2cm]
0 & -\frac{1}{2^5} & -\frac{1}{2^5} & 0 & 0 & 0\\[0.2cm]
0 & -\frac{1}{2^5} & -\frac{1}{2^5} & 0 & 0 & \frac{7}{2^{9}}\\[0.2cm]
0 & -\frac{1}{2^5} & -\frac{1}{2^5} & 0 & 0 & \frac{7}{2^{9}}\\[0.2cm]
0 & -1 & -1 & 1 & 1 & \frac{6}{2^{4}}\\[0.2cm]
\end{array}
\right),
\end{equation}
\[
X=(x_{1}, x_{2}, ... ,x_{6}).
\]
By solving it, we have 
\begin{align*}
 x_{1}=x_{6}=0,\ \ \ \ x_{3}=-x_{2},\ \ \ \  x_{5}=-x_{4}.
\end{align*}
So $\dim E^{e}(\frac{1}{16})=2$ and if we let $x_{2}=1, x_{4}=0$ or $x_{2}=0, x_{4}=1$ we can obtain a basis of $E^{e}(\frac{1}{16})$:
\begin{align*}
 b^{1}=e^{\tau_{f}}- e^{\tau_{f}\tau_{e}},\ \ \ \ b^{2}=f-f^{\tau_{e}}.
\end{align*}
Besides, By direct computation, 
we have 
\begin{align*}
 & u\cdot b^{1}=\frac{39}{2^7}b^{1}-\frac{1}{2^4}b^{2},\ \ \ \ u\cdot b^{2}=-\frac{1}{2^4}b^{1}+\frac{271}{2^7}b^{2},\\
 & v\cdot b^{1}=\frac{929}{2^9}b^{1}-\frac{1}{2^4}b^{2},\ \ \ \ v\cdot b^{2}=\frac{1}{2^4}b^{1}+\frac{1}{2^9}b^{2}.
\end{align*}
Let
\begin{align*}
 \widetilde{b^{1}}=b^{1}+l_{1}b^{2},\ \ \ \ \
 \widetilde{b^{2}}=b^{1}+l_{2}b^{2},
\end{align*}
where
\begin{align*}
 l_{1}=\frac{-29+13\sqrt{5}}{2},\ \ \ \ \
  l_{2}=\frac{-29-13\sqrt{5}}{2}.
\end{align*}
Then it is a routine work to check that 
\begin{align*}
&  \widetilde{u}\cdot\widetilde{b^{1}}=\frac{57}{32}\widetilde{b^{1}},\ \ \ \ \widetilde{u}\cdot\widetilde{b^{2}}=\frac{5}{32}\widetilde{b^{2}},\\
&  \widetilde{v}\cdot\widetilde{b^{1}}=\frac{5}{32}\widetilde{b^{1}},\ \ \ \ \widetilde{v}\cdot\widetilde{b^{2}}=\frac{57}{32}\widetilde{b^{2}}.
\end{align*}
A similar argument as Claim 1 in theorem $\ref{C3structure}$ together with the fact $\dim E^{e}(\frac{1}{16})=2$ tells us
$n_{9}=n_{10}=1$.

\emph{Claim 4) All other $n_{i}$ are equal to one.}

By direct computation, we have
\begin{align*}
&  q\dim_{[0, 0, 0]}([0, 0, 0])=q\dim_{[0, 0, 0]}([0, \frac{15}{2}, \frac{15}{2}])=q\dim_{[0, 0, 0]}([\frac{1}{2}, 0, \frac{15}{2}])
=q\dim_{[0, 0, 0]}([\frac{1}{2}, \frac{15}{2}, 0])=1,\\
&  q\dim_{[0, 0, 0]}([0, \frac{3}{4}, \frac{13}{4}])=q\dim_{[0, 0, 0]}([0, \frac{13}{4}, \frac{3}{4}])=
q\dim_{[0, 0, 0]}([\frac{1}{2}, \frac{3}{4}, \frac{3}{4}])=q\dim_{[0, 0, 0]}([\frac{1}{2}, \frac{13}{4}, \frac{13}{4}])=3+2\sqrt{2},\\
& q\dim_{[0, 0, 0]}([\frac{1}{16}, \frac{5}{32}, \frac{57}{32}])=q\dim_{[0, 0, 0]}([\frac{1}{16}, \frac{57}{32}, \frac{5}{32}])
=q\dim_{[0, 0, 0]}([\frac{1}{16}, \frac{57}{32}, \frac{165}{32}])=q\dim_{[0, 0, 0]}([\frac{1}{16}, \frac{165}{32}, \frac{57}{32}])\\
& =4+2\sqrt{2}.
\end{align*}
Let
\begin{align*}
& U^{1}=[0, 0, 0]\oplus n_{2}[0, \frac{15}{2}, \frac{15}{2}]\oplus n_{3}[0, \frac{3}{4}, \frac{13}{4}]\oplus n_{4}[0, \frac{13}{4}, \frac{3}{4}],\\
& U^{2}=n_{5}[\frac{1}{2}, 0, \frac{15}{2}]\oplus n_{6}[\frac{1}{2}, \frac{15}{2}, 0]\oplus [\frac{1}{2}, \frac{3}{4}, \frac{3}{4}]\oplus n_{8}[\frac{1}{2}, \frac{13}{4}, \frac{13}{4}],\\
& U^{3}=[\frac{1}{16}, \frac{5}{32}, \frac{57}{32}]\oplus [\frac{1}{16}, \frac{57}{32}, \frac{5}{32}]\oplus n_{11}[\frac{1}{16}, \frac{57}{32}, \frac{165}{32}] \oplus n_{12}[\frac{1}{16}, \frac{165}{32}, \frac{57}{32}].
\end{align*}
Then by a similar argument as Claim 3 and 5 in theorem $\ref{C3structure}$ we have the following results:
\begin{align*}
 q\dim_{[0, 0, 0]}U^{1}=q\dim_{[0, 0, 0]}U^{2},\ \ \ \ \ 
q\dim_{[0, 0, 0]}(U^{1}+U^{2})=q\dim_{[0, 0, 0]}U^{3},
\end{align*}
which gives the following two equations:
\begin{align}
 & 1+n_{2}+(n_{3}+n_{4})(3+2\sqrt{2})=n_{5}+n_{6}+(1+n_{8})(3+2\sqrt{2}),\label{A5qdim1}\\ 
& 1+n_{2}+n_{5}+n_{6}+(n_{3}+n_{4}+1+n_{8})(3+2\sqrt{2})=(2+n_{11}+n_{12})(4+2\sqrt{2}). \label{A5qdim2}
\end{align}
From the quantum dimension we know $[0, \frac{15}{2}, \frac{15}{2}], [\frac{1}{2}, 0, \frac{15}{2}], [\frac{1}{2}, \frac{15}{2}, 0]$ are simple current modules of $[0, 0, 0]$. By proposition 5.1 (1) in $\cite{DMZ}$ we know that $n_{2}, n_{5}, n_{6}$ are $0$ or $1$. Now assume $n_{2}=0$. 
By ($\ref{A5qdim1}$), we have 
\begin{align*}
n_{5}+n_{6}=1+n_{2}=1.
\end{align*}
So $n_{5}=1, n_{6}=0$ or $n_{5}=0, n_{6}=1$. For the first case, since 
\begin{align*}
 [\frac{1}{2}, 0, \frac{15}{2}]\boxtimes [\frac{1}{16}, \frac{57}{32}, \frac{5}{32}]=[\frac{1}{16}, \frac{57}{32}, \frac{165}{32}]
\end{align*}
A similar argument as claim 4 in theorem $\ref{C3structure}$ gives $n_{11}=1$. Then from ($\ref{A5qdim2}$) we have 
\begin{align*}
& n_{3}+n_{4}+1+n_{8}=3+n_{12},\\
& 2+3(n_{3}+n_{4}+1+n_{8})=4(3+n_{12}).
\end{align*}
From the two equations we can deduce that $2=3+n_{12}$, which is impossible! Similarly we can get a contradiction from the case $n_{5}=0, n_{6}=1$. So our assumption is wrong, which implies $n_{2}=1$. Then by ($\ref{A5qdim1}$) we can immediately get $n_{5}=n_{6}=1$. Finally, by the fusion rules
\begin{align*}
& [\frac{1}{2}, 0, \frac{15}{2}]\boxtimes [\frac{1}{2}, \frac{3}{4}, \frac{3}{4}]=[0, \frac{3}{4}, \frac{13}{4}],\\
& [\frac{1}{2}, \frac{15}{2}, 0]\boxtimes [\frac{1}{2}, \frac{3}{4}, \frac{3}{4}]=[0, \frac{13}{4}, \frac{3}{4}],\\
& [\frac{1}{2}, \frac{15}{2}, \frac{15}{2}]\boxtimes [\frac{1}{2}, \frac{3}{4}, \frac{3}{4}]=[0, \frac{13}{4}, \frac{13}{4}],\\
& [\frac{1}{2}, 0, \frac{15}{2}]\boxtimes [\frac{1}{16}, \frac{57}{32}, \frac{5}{32}]=[\frac{1}{16}, \frac{57}{32}, \frac{165}{32}],\\
& [\frac{1}{2}, \frac{15}{2}, 0]\boxtimes [\frac{1}{16}, \frac{5}{32}, \frac{57}{32}]=[\frac{1}{16}, \frac{165}{32}, \frac{57}{32}],
\end{align*}
we can obtain $n_{3}=n_{4}=n_{8}=n_{11}=n_{12}=1$.
\end{proof}
\begin{theorem}
When $\langle e, f\rangle=\langle e, e^{\tau_{f}}\rangle=\frac{3}{2^9}$, $\mathrm{VA}(e,f)$ has a unique VOA structure.
\end{theorem}
\begin{proof}
The theorem follows directly from theorem 3.6 in $\cite{DZ}$ and theorem $\ref{A5structure}$ above.
\end{proof}
\section*{Acknowledgments}
The author Wen Zheng thanks Chongying, Dong for his discussions and helpful comments.

\end{document}